\documentclass[reqno,12pt]{amsart}
\usepackage{amsmath}
\usepackage{amssymb}
\usepackage{xypic}
\usepackage{mathrsfs}
\usepackage{hyperref}
\usepackage{graphicx}
\usepackage{perpage}
\MakePerPage{footnote}
\usepackage[top=1.3in,bottom=1.3in,left=1.35in,right=1.35in]{geometry}

\date{\today}
          
 \usepackage[usenames,dvipsnames]{pstricks}
\usepackage[shell]{pdftricks}
 \begin{psinputs}
 \usepackage{pst-grad} 
 \usepackage{pst-plot} 
 \usepackage{epsfig}
 \end{psinputs}

\title[Stratifications for the surfaces $\boldsymbol{y^a=z^bx^c+x^d}$]{Whitney, Kuo-Verdier and Lipschitz stratifications for the surfaces $\boldsymbol{y^a=z^bx^c+x^d}$}

\author{Dwi Juniati, Laurent Noirel and David Trotman \vspace{-0.8cm}} 
\address{D. Juniati: Department of Mathematics, Universitas Negeri Surabaya (UNESA), Surabaya, Indonesia.}
\address{L. Noirel, D. Trotman: Aix Marseille Univ, CNRS, Centrale Marseille, I2M, Marseille, France.}

\begin{document}
\maketitle
\newtheorem{thm}{Theorem}[section]
\newtheorem{lem}[thm]{Lemma}

\theoremstyle{theorem}
\newtheorem{prop}[thm]{Proposition}
\newtheorem{defn}[thm]{Definition}
\newtheorem{cor}[thm]{Corollary}
\newtheorem{example}[thm]{Example}
\newtheorem{xca}[thm]{Exercise}

\theoremstyle{remark}
\newtheorem{rem}[thm]{Remark}

\numberwithin{equation}{section}

\newenvironment{pf}{\noindent \textbf{Proof:}}{\hfill $\square$\\}
\newcommand{\bb}{\mathbb}
\newcommand{\al}{\mathcal}
\newcommand{\ak}{\mathfrak}
\newcommand{\fs}{\mathscr}

\renewcommand{\thefootnote}{\fnsymbol{footnote}}
\begin{abstract}
We specify the canonical stratifications satisfying respectively Whitney $(a)$-regularity,  Whitney $(b)$-regularity,  Kuo-Verdier $(w)$-regularity,  and Mostowski $(L)$-regularity  for the
family of surfaces $y^a=z^bx^c+x^d$, where $a, b, c, d$ are
positive integers, in both the real and complex cases.
\end{abstract}

\parskip .12cm

\section{Introduction.}
In this paper we give diagrams showing when the regularity conditions imposed by Whitney on stratifications known as  $(a)$ and $(b)$, the condition $(w)$ introduced by Kuo and Verdier and the condition $(L)$ of T. Mostowski,  respectively hold for the stratification
with two strata given by ($V - Oz, Oz$) where
$V$ is the germ of the algebraic surface 
 $$\{(x, y, z) : y^a= z^bx^c+ x^d \} $$
in $\bb R^3$ or $\bb C^3$ for positive integers $a, b, c, d$.

The third author was present at a workshop in G\" ottingen in October 1973, when C. T. C. Wall proposed that one make such  a systematic study of the canonical $(a)$-regular and $(b)$-regular stratifications for this particular  family of surfaces. The third author then began such a study in his Warwick thesis \cite{Tr1}, finally completing the calculations required in  \cite{Tr3}.

T.-C. Kuo has  often  said that   it is important to  understand the regularity properties of the stratification of this family, in particular for constructing examples illustrating general phenomena. He gave in 1971 the first such calculations  \cite{Kuo} for the cases $(a,b,c,d) = (2k, 2, 2(k + \ell ) + 1, 2(k + \ell + m) + 1)$, with $k > 0, \ell \geq 0, m \geq 0$, in his paper introducing the ratio test $(r)$ as a sufficient criterion for Whitney $(b)$-regularity for semi-analytic strata. (Kuo's proof that $(r)$ implies $(b)$ works also for sub-analytic strata \cite{Verdier}, and in fact $(r)$ implies $(b)$ for definable stratifications in any o-minimal structure, as  shown recently by G. Valette and the third author \cite{T-V}.)

A  study of the family $\{(x, y, z) : y^a= z^bx^c+ x^d \} $, for the $(w)$-regularity of J.-L. Verdier,  was undertaken by the second author in his thesis \cite{N}, completing earlier unpublished calculations by Kuo and the third author.  Another systematic study was carried out by the first author in her thesis \cite{Juniati} for Mostowski's $(L)$-regularity. Both of these theses were directed by the third author at the University of Provence in Marseille.

 Diagrams 1-4 in section 3 below were given by the third author in \cite{Tr3} in a  publication of the University of Paris 7, with limited diffusion. Diagrams 5-6 below are from the second author's thesis \cite{N}. Diagrams 7-8 below are from  the first author's thesis \cite{Juniati}. It is clearly useful  to have the results of these systematic calculations available, and in the same place. For example using the diagrams  it becomes straightforward to specify infinitely many real algebraic examples where $(b)$-regularity holds but $(w)$ fails. Previously only three such real algebraic examples had been found, in the 1979 paper of H. Brodersen and the third author \cite{BT}. These correspond to the choices $(a,b,c,d) \in \{ (3,2,3,5), (4,4,1,3), (4,2,5,7) \}.$ Note that Verdier explicitly stated in early  1976 \cite{Verdier} that he knew of no sub-analytic example satisfying $(b)$ but not $(w)$. The first semi-algebraic example was constructed at Oslo in August 1976 by the third author \cite{Tr2}.
 
  The continued need for more available information about real algebraic examples, as provided here, is illustrated by a recent paper by Cluckers, Comte and Merle \cite{CCM}, where the authors mention that $(b)$ does not imply $(w)$ for real algebraic varieties, but give no reference. Ren\'e Thom  proposed around 1985 that more needed to be understood concerning the difference between $(b)$ and $(w)$. 
  
  Another area in which it is important to understand the difference between Whitney $(b)$-regularity and Kuo-Verdier $(w)$-regularity is the study of continuity of the density along strata of a definable regular stratification. For definable sets in arbitrary o-minimal structures $(w)$-regularity implies continuity of the density (in fact the density is a Lipschitz function as shown in  \cite{Va}  and \cite{NV2}) while the density can fail to be continuous along a stratum of a $(b)$-regular stratification of a definable set if the o-minimal structure is not polynomially bounded \cite{T-V}. Yet another difference was shown by Navarro Aznar and the third author for subanalytic stratifications: $(w) \implies (w^*)$, but $(b) \implies (b^*)$ is unknown except when the small stratum has dimension one \cite{NT}.
 
\section{Definitions.}

The notion of stratification is one of the most fundamental concepts in
algebraic geometry and  singularity  theory. The idea of decomposing a
singular topological space $X$ into non singular parts which we call
the strata, satisfying some regularity conditions, goes back to
Whitney and Thom's work. They suggested the use of stratifications
 as a method of understanding the geometric structure of singular
analytic spaces. 

At first, Whitney defined the stratification of a real algebraic set
by iteratively taking the singular parts. The problem with this
stratification is the absence of topological triviality along
strata. Later Whitney introduced two famous regularity conditions for a
stratification,  called Whitney $(a)$ and $(b)$-regularity.

These regularity conditions are defined in the following way :

\begin{defn} Let $M$ be a smooth manifold, and let $X, Y$ be disjoint smooth
submanifolds of $M$ such that $ Y \cap \overline {X} \neq \emptyset$.

 (i). The pair $(X,Y)$ is said to be $(a)$-regular at a point $y \in Y \cap \overline{X}$ if for each sequence  of points $\{x_i\}$ in X converging to $y$ such that the sequence of tangent spaces $\{T_{x_i}X \}$ converges to
$\alpha$, then $T_yY$ $\subset \alpha$.

(ii) The pair $(X,Y)$ is said to be $(b)$-regular at a point $y \in Y \cap \overline{X}$ if for each pair of sequences of points $\{x_i\}$ in X and $\{y_i\}$ in Y converging to $y$
such that the sequence of tangent spaces $\{T_{x_i}X \}$ converges to
$\alpha$, and the sequence of unit vectors $\overline {x_i y_i}$ converges to $\lambda$, then $\lambda$ $\subset \alpha $.
\end{defn}

Whitney showed that all analytic varieties in $\bb C^n$ have a stratification satisfying
these conditions. Thom used the Whitney conditions in the differentiable case and proved that $(b)$-regular stratifications are topologically trivial along strata; this is  called 
Thom's isotopy Theorem. 

Another important condition is the $(w)$-regularity introduced by T.-C. Kuo \cite{Kuo} and J.-L. Verdier \cite{Verdier}.

\begin{defn}
The pair $(X,Y)$ as above is said to be $(w)$-regular at $y_0 \in Y \cap \overline{X}$ if there is a real constant $C>0$ and a neighborhood $U$ of $y_0$ such that
$$\delta (T_yY,T_xX) \leq  C |x-y|$$
 for all  $y \in U \cap Y$ and all  $x \in U \cap X.$
 \end{defn}
 
 Using Hironaka's resolution of singularities, Verdier established
that every analytic variety, or subanalytic set, admits a $(w)$-regular stratification. He also showed that $(w)$-regularity of a subanalytic stratification implies $(b)$-regularity. This proof was then generalised by T\`a L\^e Loi to the general case of definable stratifications in an arbitrary o-minimal structure \cite{Loi}. In
general, for $C^2$ differentiable stratifications,  $(w)$-regularity also
implies local topological triviality along strata, see \cite{Verdier}.

In his thesis \cite{Tr1}, the third author gave the first proof that $(b)$-regularity is strictly
weaker than $(w)$-regularity for both semi-algebraic stratified sets and real algebraic sets. But $(w)$-regularity is not strong enough to give  the stability of such
characteristics of singularities as order of contact. In \cite{P4},
Parusinski gave an example of a stratification that satisfies $(w)$-regularity but where the order of contact is not preserved. Another example, due to S. Koike, was described in \cite{Juniati3}.

The stronger the conditions imposed on the stratified set, the
better the understanding of the geometry of this set. One of the
strongest generic conditions is the $(L)$-regularity (defined below)  introduced by T.
Mostowski in 1985 \cite{Mostowski}. Mostowski introduced this notion of Lipschitz
stratification and proved the existence of such stratifications for complex analytic sets.
The existence of Lipschitz stratification for real analytic sets, then for semi-analytic sets and subanalytic sets, was later proved by A. Parusinski in \cite{P1, P2, P3}. Recently N. Nguyen and G. Valette have proved the existence of Lipschitz stratifications of definable sets in polynomiallyt bounded o-minimal structures \cite{NV}.  Lipschitz stratifications ensure local bi-Lipschitz triviality
of the stratified set along each stratum, and bi-Lipschitz
homeomorphisms preserve sets of measure zero, order of contact, and
Lojasiewicz exponents. The condition $(L)$ of Mostowski is preserved after
intersection with generic wings, that is $(L)$-regularity implies
$(L^*)$-regularity, see \cite{Juniati3}; this was one of the criteria imposed on any good equisingularity notion by B. Teissier in his foundational
1974 Arcata paper \cite{Te}. We recall now the definition of Lipschitz
stratification due to Mostowski, where we simplify slightly the (equivalent) formulation of the definition.

 Let $X$ be a closed subanalytic subset of an
 open subset of $\bb R^n$.
By a {\it stratification} of $X$ we shall mean a family $\Sigma = \{
S^j\}_{j=l}^m$ of closed subanalytic subsets of $X$ defining a
filtration:
$$X = S^m\supset S^{m-1}\supset \cdots \supset S^l \neq  \emptyset$$
and  ${\mathring {S^j}} = S^j \setminus S^{j-1}$, for  $j=l,\,l+1,\dots,\,m $ (where
$S^{l-1}=\emptyset$),
 is a smooth manifold of pure dimension $j$ or empty.
We call the connected components of
 $\mathring{S^j}$ the {\it strata} of $\Sigma$.
We denote the function measuring distance to $S^j$
 by $ d_j$, so that $d_j(q)$ = dist$(q,S^j)$. Set $d_{l-1}\equiv 1$, by convention
(this will be used in the definitions below). 

\begin{defn}
Let $\gamma >1$ be a fixed constant. A {\it chain} for a point $q \in \mathring{S^j}$  is a (strictly) decreasing sequence of indices $j = j_1,  j_2, \cdots, j_r=l$ such that each $j_s$ for $s \geq 2$ is the largest integer less than
$j_{s-1}$ for which
$$ d_{{j_s}-1}(q) \geq 2\gamma^2 d_{j_s}(q) .$$
For each $j_s \in \{j_1, \dots, j_r \}$ choose a point
 $q_{j_s} \in {\mathring{S^{j_s}}}$ such that  $q_{j_1}=q$  and $|q-q_{j_s}| \leq \gamma d_{j_s}(q). $
\end{defn}
If there is no confusion, we will call the sequence of
 points $(q_{j_s})$  a chain of $q$.
 
 For $q \in \mathring{S^j}$,  let $P_q: \bb R^n \rightarrow T_q \mathring{S^j}$ be
 the orthogonal projection to the tangent space $T_q \mathring{S^j}$ and let
$P_q^{\perp} = I-P_q$ be the orthogonal projection onto the normal
space  $(T_q\mathring{S^j})^{\perp}$.

A stratification $\Sigma =\{ S^j\}^m_{j=l}$ of $X$ is said to be a {\it Lipschitz stratification} if for some constant $C>0$ and for
every chain $q=q_{j_1}, q_{j_2}, \cdots, q_{j_r}$ and every
 $k$ for $2\leq k \leq r$,
\begin{equation}
|P^{\perp}_{q}P_{q_{j_2}}\cdots P_{q_{j_k}}| \leq C|q-q_{j_2}| / d_{{j_k}-1}(q) \tag{L1}
 \end{equation}
and for each $q'\in \mathring{S^{j_1}}$ such that $|q-q'| \leq(1/2\gamma) d_{{j_1}-1}(q)$,
\begin{equation}
|(P_q -P_{q'})P_{q_{j_2}}\cdots P_{q_{j_k}}| \leq C |q-q'| /d_{{j_k}-1}(q) \tag{L2}
\end{equation}
and
\begin{equation}
|P_q - P_{q'}| \leq C |q-q'| / d_{{j_1}-1}(q) \tag{L3}
\end{equation}

\medskip

T. Mostowski and A. Parusi\'nski  showed that Lipschitz stratifications satisfy a bilipschitz version  of  the Verdier isotopy theorem.

It is easy to see that condition (L1)  in the definition of Lipschitz stratification of Mostowski implies 
$(w)$-regularity and that in our case, where there are only two strata, it is
actually equivalent to $(w)$-regularity.

We observe that the algebraic surfaces  $\{ y^a = t^b x^c + x^d \}$ are quasi-homogeneous with weights $(a, d, \frac{bd}{d-c})$, when $c < d$, which is the situation where the different types of regularity are in doubt.
 In stratification theory, one often studies quasi-homogeneous spaces. S.
Koike in \cite{Koike1} and \cite{Koike2}, studied modified Nash
triviality in  the case of quasi-homogeneous families with
isolated singularities.
Fukui and Paunescu  \cite{F-P}  found  some  conditions on the
equations defining the strata to show that a given family is
topologically trivial. They gave a weighted version  of  $(w)$-regularity and of  Kuo's ratio test, thus providing weaker sufficient conditions
for topological stability. Sun \cite{Sun} used these methods to generalize the work
 of the third author  and Wilson \cite{T-W} for $(t)$-regularity.
The first author and G. Valette  in \cite{Juniati4} gave some conditions for a
quasi-homogeneous stratification to satisfy $(w)$-regularity or
to be bi-Lipschitz trivial (which means to satisfy the conclusion of
the bi-Lipschitz Isotopy Lemma). They  also studied the
variation of the volume of quasi-homogeneous stratified families.

 \section {\bf{ The classifications.} }

In this section we give diagrams showing when $(E)$-regularity, with $E$ successively $a$, $b$, $w$ and $L$,  holds for the stratification
with two strata given by ($V - Oz, Oz$) where
$V$ is the germ of the algebraic surface 
 $$\{(x, y, z) : y^a= z^bx^c+ x^d \} $$
in $\bb R^3$ or $\bb C^3$ for positive integers $a, b, c, d$.

\subsection{Classification  for Whitney (a)-regular and (b)-regular stratifications}

The long and detailed calculations for these stratifications can be
seen in chapter 8 of the third author's thesis \cite{Tr1}. He gave
the outstanding calculations for these in \cite{Tr3}.

Whitney $(a)$-regularity  is verified at $(0,t_0)$ in $\bb R^n
\times \bb R^k$ for $(X,Y)$ when 
$$ \frac{\frac{\partial{f}}{\partial{t_{i}}}(x,t)}{|{\frac{\partial{f}}{\partial{x_{1}}} ,\ldots, \frac{\partial{f}}{\partial{x_n}} |}}$$
tends to $0$ as $(x,t)$ tends to $(0,t_0)$ on $X$ for $1 \leq i \leq k.$

For Whitney $(b)$-regularity, Trotman used $(b^{\pi})$-regularity. The pair $(X,Y)$ is said to be $(b^{\pi})$-regular $y \in Y \cap \overline{X}$ when every limit of
tangent spaces to $X$ at $y$ defined by a sequence $\{ x_i \}$ in $X$
contains the limiting direction of secants $x_i \pi (x_i) $ where
$\pi$ is a (linear) retraction onto $Y$ of a tube around $Y$ defined
near $y$. When both $(a)$ and $(b^{\pi})$ hold, $X$ is $(b)$-regular. And for the calculation, th
$(b^{\pi})$-condition is verified at $(0,t_0)$ in ${\Bbb {R}}^n
\times {\Bbb {R}}^k$ for $(X,Y)$ when 
$$\frac{\Sigma_{j=1}^{n}  x_j \frac{\partial{f}}{\partial{x_{j}}}(x,t)}   {| (x_1,\ldots, x_n)| . |
({\frac{\partial{f}}{\partial{x_{1}}}, \ldots , \frac
{\partial{f}}{\partial{x_n}}, \frac{\partial{f}}{\partial{t_1}},
..., \frac{\partial{f}}{\partial{t_k}} )| }}$$
tends to $0$ as
$(x,t)$ tends to $(0,t_0)$ in $X$.

 It is easily verified using the curve selection lemma for
semi-algebraic sets that this condition holds if and only if  the
limits tend to $0$ for $(x,t)$ lying on an analytic curve $\alpha (u) =
(x(u),t(u))$ such that $\alpha (0) = (0,t_0)$ and $\alpha (u) \in X$
if $u \neq 0$. The third author used this fact in some of the calculations
required to complete the diagrams, but he used also a                                                                                                                                                                                                                                                                                                                                                                                                                                                                                                                                                                                                                                                                                                                                                                                                                                                                                                                                                                                                                                                                                                                                                                                                                                                                                                                                                                                                  simpler approach by considering separately sequences $(x_i,t_i)$ whose
limiting direction is either tangent to $0 \times \bb R^k$,
tangent to $\bb R^n \times 0$ or neither of these, making three
separate cases in all.

He distinguished between three distinct cases considering 
sequences $(x_i, t_i)$ tending to $(0,0)$ with 

(i) $\lim_{i \to \infty} \{\frac {x_i}{t_i}\} = 0$.

(ii) $ \lim_{i \to \infty} \{ \frac {x_i}{t_i}\} = \kappa \neq 0$, $ \kappa \in \bb R$.

(iii) $ \lim_{i \to \infty} \{ \frac {t_i}{x_i}\} = 0$.

 The difference between the diagrams in the real and complex cases arises solely from
the non-existence of branches of the curve $\{z^b +x^{d-c}=0\}$ near
0 in $\bb R^2$ when $b$ and $d-c$ are even and $d>c$. Such
branches exist in $\bb C^2$, and they also
exist for the curve $\{-z^b +x^{d-c}=0\}$ in $\bb R^2$.

\subsection{Classification for Kuo-Verdier (w)-regular stratifications}

It is known that $w$-regularity implies 
Whitney $b$-regularity \cite{Kuo}, \cite{Verdier}. So we only need to study cases when $b$-regularity is satisfied in the classification obtained by Trotman.

 If $v=(0,0,1)$ then $\delta (Y,T_{(x,y,z)}X)= \frac {|\Delta
 f(x,y,z)v |}{|| \Delta f(x,y,z)||}$. Let $\pi_2$
 be the projection defined by the second coordinate. For $(x,z) \in \pi_2
 (X)$, let $Z(x,z)=\frac {|\Delta
 f(x,y,z)v |}{|| \Delta f(x,y,z)|| . || (x,y) ||}$ where $y$ is such that $(x,y,z) \in V$ with $V$ the graph of the smooth map $(x,z) \mapsto z^b x^c + x ^d $.
We show that $Z$ is independent from the determination of $y$. So the
$(w)$-regularity is verified at $0$ if $Z$ is bounded in the neighbourhood
of $0$ in $\pi_2(X)$. This means that $Z(x,z)$ converges to $0$ when
$(x,t) \in \pi_2(X)$ converges to $0$. The calculation uses the fact that $ Z \equiv Z' $, where:
$$ Z'(x,z) = \frac {| z^{b-1}x^c |} { \sup ( | cx^{c-1}z^b + d x^{d-1} |, | z^b x^c + x^d |^{\frac {a-1}{a}} ) . \sup ( | x |, | z^b x^c + x^d |^{\frac {1}{a}}
)}$$

For the detailed calculations see the second author's thesis [N].

\subsection{Classification for Lipschitz (L)-regular stratifications}
Since the condition (L1) in the definition of
Lipschitz stratification of Mostowski implies the $w$-regularity, and in our case (where there only two strata) is
 actually equivalent,  we need only to study cases when $w$-regularity is satisfied in the
 classification obtained by the second author \cite{N}. And for these we  need to verify the conditions (L2) and (L3) for
 $k = 2$.
 
 The essential technique in verifying (L2) and (L3) is
 to apply the mean value theorem to
compare values at two points $q$ and $q'$ whose distance apart
is controlled. To verify  (L2), we must prove:

\begin{align*}
 |(P_{q } - P_{q'})
P_{q_{j_2}}(v) | \,
&=  |(P_{q} - P_{q'})(0,0,1) | \\
&= | (P_{q}^\perp -P_{q'}^\perp) (0,0,1) \\
&= \left |{\frac { bx^c\, z^{b-1} \partial_{q} f} {||
\partial_{q}f|| ^2} - \frac{ b{x'}^c\, {z'}^{b-1}
\partial_{q'}f}{\|| \partial_{q'} f||^2} } \right |\\
&\leq \, C\, |q - q'|
\end{align*}

To verify the condition (L3), we use the vector basis $v_1 = ( 1, 0, 0 ), v_2 = ( 0, 1, 0 )$ and $ v_3 = ( 0, 0, 1 )$. We prove that :
\begin{align*}
  |(P_{q}-P_{q'}) (v_1)| 
  &= \left | \frac {\frac {\partial \varphi}{\partial x}(\frac {\partial
\varphi}{\partial x}, \frac {\partial \varphi}{\partial y },
\frac {\partial \varphi}{\partial z})} { {|| \partial_q
\varphi
||} ^2} - \frac {\frac {\partial \varphi}{\partial {x'}} (
\frac {\partial \varphi}{\partial {x'}}, \frac {\partial
\varphi}{\partial {y'} }, \frac {\partial \varphi}{\partial
{z'}})} { {|| \partial_{q'} \varphi ||} ^2} \right |\\
& \leq \, C \, \frac {| q-q' |}{d(q, O_z)},
\end{align*}

\begin{align*}
  | (P_{q}-P_{q'}) (v_2) | 
  &=  \left | 
 \frac {\frac {\partial \varphi}{\partial y} (\frac {\partial
\varphi}{\partial x}, \frac {\partial \varphi}{\partial y },
\frac {\partial \varphi}{\partial z})} { {|| \partial_q
\varphi
||} ^2}- \frac {\frac {\partial \varphi}{\partial {y'}} (
\, \frac {\partial \varphi}{\partial {x'}},\, \frac {\partial
\varphi}{\partial {y'} }, \frac {\partial \varphi}{\partial
{z'}})} { {|| \partial_{q'} \varphi ||} ^2} \right
|\\
& \leq \, C \, \frac {| q-q' |}{d(q, O_z)}, 
\end{align*}
and
\begin{align*}
  | (P_{q}-P_{q'}) (v_3) |
  &= \left | 
 \frac {\frac {\partial \varphi}{\partial z} ( \frac {\partial
\varphi}{\partial x}, \frac {\partial \varphi}{\partial y },
\frac {\partial \varphi}{\partial z})} { {|| \partial_q
\varphi
||} ^2}  -  \frac {\frac {\partial \varphi}{\partial {z'}} ( \frac {\partial \varphi}{\partial {x'}},\, \frac {\partial
\varphi}{\partial {y'} }, \frac {\partial \varphi}{\partial
{z'}})} { {|| \partial_{q'} \varphi ||} ^2}
\right |\\
& \leq  C  \frac {|q-q'|}{d(q, O_z)}.
\end{align*}

To check that the condition $(L)$ fails, we give the two curves of $q$
and $q'$ where the condition $L2 $ or $L3 $ fails.
Long and detailed calculations deciding several branches of the Diagram
(see below) were given in the first author's thesis \cite{Juniati}.

\newpage

\begin{center}

\section{\bf The diagrams of the classification for Whitney, Kuo-Verdier and Lipschitz stratifications}
\end{center}
\begin{align*}
\text { Note : } & \surd \text{means that regularity holds }\\
& \times \text {means that regularity fails } \\
& ? \text { means undecided }\\
 &  b \equiv 0(2)  \text { means that b
is even }\\ & b \equiv 1(2) \text {means that b is odd, etc. }
\end{align*}

\begin{equation*}
\begin{cases}
a=1 \,\,\, \surd \\ \\ 
a>1
\begin{cases}
d \leq c \,\,\, \surd \\ \\
d > c
\begin{cases}d \geq b+c \,\,\, \times\\ \\
d < b+c
\begin{cases}
a \leq b \,\,\, \surd\\ \\
a > b \begin{cases}
d < \frac{ac}{a-b} \,\,\, \surd \\ \\
d \geq \frac{ac}{a-b}\,\,\, \times
\end{cases}
\end{cases}
\end{cases}
\end{cases}
\end{cases}
\end{equation*}

\vskip1cm

\centerline {{\bf Diagram 1}. Whitney $(a)-$regularity in ${\Bbb{C}}^3 $}

\begin{equation*}
\begin{cases}
a = 1 \,\,\, \surd \\ \\
a>1 \begin{cases}
d \leq c \,\,\, \surd \\ \\
c < d<b+c \begin{cases}
a \leq b \,\,\, \surd \\ \\
b < a
\begin{cases}
d < \frac{ac}{a-b} \,\,\, \surd \\ \\
d \geq \frac{ac}{a-b} 
\begin{cases}
b \equiv 0(2) \begin{cases}
d \equiv c(2) \,\,\, \surd \\ \\
d \equiv c + 1(2) \,\,\, \times
\end{cases}\\
b\equiv 1(2) \,\,\, \times \\
\end{cases}
\end{cases}
\end{cases}\\
b+c \leq d
\begin{cases}
b \equiv 1(2) \,\,\, \times \\ \\
b\equiv 0(2)
\begin{cases}
d \equiv c + 1(2) \,\,\, \times \\ \\
d \equiv c(2)
\begin{cases}
a \leq b \,\,\, \surd \\ \\
b<a<b+c
\begin{cases}
d< \frac{ac}{a-b}\,\,\,\surd \\ \\
d \geq \frac{ac}{a-b} \,\,\, \times
\end{cases}\\
b + c \leq a \,\,\, \times
\end{cases}
\end{cases}
\end{cases}
\end{cases}
\end{cases}
\end{equation*}
\vskip 1cm
\centerline {{\bf Diagram 2}. Whitney $(a)-$regularity in ${\Bbb{R}}^3 $}

\begin{equation*}
\begin{cases}
a = 1 \,\,\, \surd \\ \\
a > 1 \begin{cases}
d \leq c \,\,\, \surd \\ \\
d > c \,\,\, \times
\end{cases}
\end{cases}
\end{equation*}
\vskip .5cm

\centerline {{\bf Diagram 3}. Whitney $(b)-$regularity in ${\Bbb{C}}^3 $ }

\vskip 2cm

\begin{equation*}
\begin{cases}
a = 1 \,\,\, \surd \\ \\
a > 1 \begin{cases}
d \leq c \,\,\, \surd \\ \\
d > c \begin{cases}
b \equiv 1(2) \,\,\, \times \\ \\
b \equiv 0(2) \begin{cases}
d \equiv c + 1(2) \,\,\, \times \\ \\
d \equiv c(2) \begin{cases}
d < a \begin{cases}
d < b+c \,\,\, \surd \\ \\
d \geq b+c \,\,\, \times
\end{cases}\\ \\
a \leq d
\begin{cases}
a > c \,\,\, \times \\ \\
a \leq c \begin{cases}
d < b+c \,\,\, \times \\ \\
d \geq b+c \begin{cases}
a \leq b \,\,\, \surd \\ \\
a > b \begin{cases}
d < \frac{ac}{a-b} \,\,\, \surd \\ \\
d \geq \frac{ac}{a-b} \,\,\, \times 
\end{cases}
\end{cases}
\end{cases}
\end{cases}
\end{cases}
\end{cases}
\end{cases}
\end{cases}
\end{cases}
\end{equation*}

\vskip .5cm

\centerline {{\bf Diagram 4}. Whitney $(b)-$regularity in ${\Bbb{R}}^3 $}


\begin{equation*}
\begin{cases}
a = 1 \,\,\, \surd \\ \\
a > 1 \begin{cases}
d \leq c \,\,\, \surd \\ \\
d > c \,\,\, \times
\end{cases}
\end{cases}
\end{equation*}

\vskip .5cm

\centerline {{\bf Diagram 5}. Kuo-Verdier $(w)-$regularity in ${\Bbb{C}}^3 $ }

\begin{equation*}
\begin{cases}
a = 1 \,\,\, \surd \\ \\
a > 1 \begin{cases}
d \leq c \,\,\, \surd \\ \\
d > c \begin{cases}
b \equiv 1(2) \,\,\, \times \\ \\
b \equiv 0(2) \begin{cases}
d \equiv c + 1(2) \,\,\, \times \\ \\
d \equiv c(2) \begin{cases}
a(d-c) \leq b | d-a| \,\,\, \surd \\ \\
a(d-c)>b|d-a| \,\,\, \times
\end{cases}
\end{cases}
\end{cases}
\end{cases}
\end{cases}
\end{equation*}

\vskip .5cm
\centerline {{\bf Diagram 6}. Kuo-Verdier $(w)-$regularity in ${\Bbb{R}}^3 $ }

\begin{equation*}
\begin{cases}
a =1 \,\,\, \surd \\ \\
a > 1 \begin{cases}
d \leq c \,\,\, \surd \\ \\
d > c \,\,\, \times
\end{cases}
\end{cases}
\end{equation*}

\vskip .5cm

\centerline {{\bf Diagram 7}. Mostowski $(L)-$regularity in ${\Bbb{C}}^3 $ }

\pagebreak
\centerline {{\bf Diagram 8}. Mostowski $(L)-$regularity in ${\Bbb{R}}^3 $ }
\begin{equation*}
\begin{cases}
a=1 \,\,\, \surd \\ \\
a > 1 \begin{cases}
d \leq c \,\,\, \surd \\ \\
d > c \begin{cases}
b \equiv 1(2) \,\,\, \times \\ \\
b \equiv 0(2) \begin{cases}
d \equiv c + 1(2) \,\,\, \times \\ \\
d \equiv c (2) \begin{cases}
a(d-c) \leq b |d-a| \begin{cases}
d \geq a \,\,\,\,(\text{continued below}) \\ \\
d<a \begin{cases}
d,c \,\,\, \text{even} \,\,\,\,\, \times \\\\
d,c \,\,\, \text{odd} \begin{cases}
b < 2 (d-c) \,\,\, \times \\ \\
b \geq 2(d-c) \,\,\, ?
\end{cases}
\end{cases}
\end{cases}\\ \\
a(d-c) > b |d-a|\,\,\, \times
\end{cases}
\end{cases}
\end{cases}
\end{cases}
\end{cases}
\end{equation*}

\begin{equation*}
\rightarrow d \geq a \begin{cases}
a \,\,\, \text{even} \,\,\, \times \\ \\
a \,\,\, \text{odd} \begin{cases}
c < a \,\,\, \times \\ \\ 
c \geq a \begin{cases}
b<a\begin{cases}
a (d-c) \leq d - a \,\,\, \surd \\ \\
a(d-c) > d-a \begin{cases}
2a(d-c) > db \,\,\, \times \\ \\ 
2a(d-c) \leq db \,\,\, ?
\end{cases}
\end{cases}\\
b\geq a \begin{cases}
b \geq 2a \,\,\, \surd \\ \\
b<2a \begin{cases}
a(d-c)\leq d \,\,\, \surd \\ \\
a(d-c) > d \begin{cases}
2a(d-c)>db \,\,\, \times \\ \\
2a(d-c)\leq db \,\,\, ?
\end{cases}
\end{cases}
\end{cases}
\end{cases}
\end{cases}
\end{cases}
\end{equation*}

\section{Some consequences.}

An immediate consequence of diagram 1, determining the canonical $(a)$-regular stratification of $V(a,b,c,d)$ in ${\Bbb{C}}^3$ is :

\medskip

{\bf Proposition.} {\it There is no algebraic invariant ${\alpha}$ of plane curves in ${\Bbb{C}}^2$ such that   ${\alpha}(F_t)$ is constant if and only if $(F^{-1}(0) - (0 \times {\Bbb{C}}), 0 \times {\Bbb{C}})$ is $(a)$-regular where $F : {\Bbb{C}}^2 \times {\Bbb{C}}, 0 \times {\Bbb{C}} \rightarrow {\Bbb{C}}, 0$ is a polynomial map whose singular set is $0 \times {\Bbb{C}}$, and $F_t(z) = F(z,t)$.}

{\bf{Proof}}. It is enough to consult diagram 1 where one sees that $(a)$-regularity depends on $b$, which governs the speed of the deformation in $J^k(2,1)$. If such an invariant existed, $(a)$-regularity would be independent of $b$.

\medskip

{\bf{Notes.}} 

1. The Milnor number is an algebraic invariant whose constance is equivalent to $(b)$-regularity. Diagram 3 confirms this : it is independent of the power $b$ of $t$ in the deformation $- y‰ + t^b x^c + x^d$.

2. By an example of Zariski \cite{Zariski}, ${ y^2 = st x^2 + x^3}$, one knew there was no semicontinuous invariant $\alpha$. But there do exist nonsemicontinuous numerical invariants of plane curves (cf. chapter 10 of Milnor's book \cite{Milnor}).

3. When $b=1$, there is an $(a)$-fault at $0$ if and only if there is a $(b)$-fault at 0. This follows from diagrams 1-4 and is rather surprising.

\medskip

{\bf{$C^1$-smooth varieties.}}

An algebraic variety can be $C^1$-smooth and yet be singular. In fact for all finite $k$ there is a singular algebraic variety which is $C^k$-smooth. However J. Mather has proved \cite{Mather} that if a variety is $C^{\infty})$-smooth, then it is nonsingular. Obviously the converse is true too. The relevance to our study follows from the $C^1$-invariance of $(a)$- and $(b)$-regularity, which implies that both $(a)$ and $(b)$ hold if $V$ is $C^1$-smooth. The following proposition, whose proof is an easy exercise, may be compared with diagrams 2 and 4. The proposition applies when $d > c$, but when $d \leq c$ we have already seen that $(a)$ and $(b)$ hold.

\medskip

{\bf{Proposition.}}  {\it Let $y‰ = t^bx^c + x^d$ and let $d > c$. Then $y$ is a $C^1$ function of $x$ and $t$ if and only if $a \leq c, d < ac/(a-b)$, $b$ is even, and $(d-c)$ is even.}

\medskip

When $f(x,t) = (t^bx^c + x^d)^{1/a}$ is a $C^1$ function of $x$ and $t$, the graph of $f$ is a $C^1$ submanifold, so that $(b)$-regularity holds.

\medskip

{\bf Examples of $(b)$-regular $(w)$-faults.}

As noted earlier  \cite{BT}  provides three examples among the surfaces $y‰ = t^b x^c + x^d$ where $(b)$-regularity holds but $(w)$-regularity fails : $(a,b,c,d) : (3,2,3,5), (4,2,5,7)$ and $(4,4,1,3)$.  The diagrams of the present paper furnish infinitely many such examples, as by Diagram 6 $(w)$ fails when $a > 1, d > c, b \equiv 0 (2), d \equiv c (2),$ and $a(d-c) > b \vert d-a \vert$, and by Diagram 4 $(b)$ holds if further $a < c, d \geq b+c $ and either $a<b $ (case 1) or $a\geq b$ and $d(a-b) < ac$ (case 2). The examples (3,2,3,5) and (4,2,5,7) correspond to case 1, while (4,4,1,3) corresponds to case 2.

\end{document}